\documentclass{amsart}

\usepackage{fullpage}
\usepackage{amsmath}
\usepackage{amsfonts}
\usepackage{amssymb}

\input xy
\xyoption{matrix}\xyoption{arrow}\xyoption{curve}  \xyoption{frame}

\begin{document}

%%%%% Functors %%%%%

\newcommand{\Hom}{\mathrm{Hom}}
\newcommand{\RHom}{\mathrm{RHom}^*}
\newcommand{\HOM}{\mathrm{HOM}}
\newcommand{\stHom}{\underline{\mathrm{Hom}}}
\newcommand{\Ext}{\mathrm{Ext}}
\newcommand{\Tor}{\mathrm{Tor}}
\newcommand{\HH}{\mathrm{HH}}
\newcommand{\Endo}{\mathrm{End}}
\newcommand{\ENDO}{\mathrm{END}}
\newcommand{\stEnd}{\mathrm{\underline{End}}}
\newcommand{\Tr}{\mathrm{Tr}}
\newcommand{\GL}{\mathrm{GL}}

%%%%% Functions/Operators (no space after) %%%

\newcommand{\coker}{\mathrm{coker}}
\newcommand{\aut}{\mathrm{Aut}}
\newcommand{\op}{\mathrm{op}}
\newcommand{\add}{\mathrm{add}}
\newcommand{\ADD}{\mathrm{ADD}}
\newcommand{\ind}{\mathrm{ind}}
\newcommand{\rad}{\mathrm{rad}}
\newcommand{\soc}{\mathrm{soc}}
\newcommand{\ann}{\mathrm{ann}}
\newcommand{\im}{\mathrm{im}}
\newcommand{\chr}{\mathrm{char}}
\newcommand{\pdim}{\mathrm{p.dim}}

%%%%% Categories %%%%%%%%%

\newcommand{\rmod}{\mbox{mod-}}
\newcommand{\Rmod}{\mbox{Mod-}}
\newcommand{\lmod}{\mbox{-mod}}
\newcommand{\lMod}{\mbox{-Mod}}
\newcommand{\stmod}{\mbox{\underline{mod}-}}
\newcommand{\stlmod}{\mbox{-\underline{mod}}}
\newcommand{\Mod}{\mathcal{M}\mbox{od}}

\newcommand{\gmod}[1]{\mbox{mod}_{#1}\mbox{-}}
\newcommand{\gMod}[1]{\mbox{Mod}_{#1}\mbox{-}}
\newcommand{\Bimod}[1]{\mathrm{Bimod}_{#1}\mbox{-}}

\newcommand{\proj}{\mbox{proj-}}
\newcommand{\lproj}{\mbox{-proj}}
\newcommand{\Proj}{\mbox{Proj-}}
\newcommand{\inj}{\mbox{inj-}}
\newcommand{\coh}{\mbox{coh-}}
\newcommand{\gr}{\mbox{gr-}}
\newcommand{\stgr}{\mbox{\underline{gr}-}}

%%%%% Style %%%%%%%%%%

\newcommand{\und}[1]{\underline{#1}}
\newcommand{\gen}[1]{\langle #1 \rangle}
\newcommand{\floor}[1]{\lfloor #1 \rfloor}
\newcommand{\ceil}[1]{\lceil #1 \rceil}
\newcommand{\bnc}[2]{\left(\scriptsize \begin{array}{c} #1 \\ #2 \end{array} \right)}
\newcommand{\bimo}[1]{{}_{#1}#1_{#1}}
\newcommand{\ses}[5]{\ensuremath{0 \rightarrow #1 \stackrel{#4}{\longrightarrow} 
#2 \stackrel{#5}{\longrightarrow} #3 \rightarrow 0}}
\newcommand{\A}{\mathcal{A}}
\newcommand{\B}{\mathcal{B}}
\newcommand{\C}{\mathcal{C}}
\newcommand{\D}{\mathcal{D}}
\newcommand{\M}{\mathcal{M}}
\newcommand{\T}{\mathcal{T}}
\newcommand{\s}{\mathcal{S}}
\newcommand{\tC}{\tilde{\mathcal{C}}}
\newcommand{\tK}{\tilde{K}(A)}
\newcommand{\ul}[1]{\underline{#1}}

%%%%% Theorems %%%%%%%%%%

\newtheorem{therm}{Theorem}[section]
\newtheorem{defin}[therm]{Definition}
\newtheorem{propos}[therm]{Proposition}
\newtheorem{lemma}[therm]{Lemma}
\newtheorem{coro}[therm]{Corollary}

\title{Tilting mutation of weakly symmetric algebras and stable equivalence}
\author{Alex Dugas}
\address{Department of Mathematics, University of the Pacific, 3601 Pacific Ave, Stockton CA 95211, USA}
\email{adugas@pacific.edu}


\begin{abstract} We consider tilting mutations of a weakly symmetric algebra at a subset of simple modules, as recently introduced by T. Aihara.  These mutations are defined as the endomorphism rings of certain tilting complexes of length $1$.  Starting from a weakly symmetric algebra $A$, presented by a quiver with relations, we give a detailed description of the quiver and relations of the algebra obtained by mutating at a single loopless vertex of the quiver of $A$.  In this form the mutation procedure appears similar to, although significantly more complicated than, the mutation procedure of Derksen, Weyman and Zelevinsky for quivers with potentials.  By definition, weakly symmetric algebras connected by a sequence of tilting mutations are derived equivalent, and hence stably equivalent.  The second aim of this article is to study these stable equivalences via a result of Okuyama describing the images of the simple modules.   As an application we answer a question of Asashiba on the derived Picard groups of a class of self-injective algebras of finite representation type.  We conclude by introducing a mutation procedure for maximal systems of orthogonal bricks in a triangulated category, which is motivated by the effect that a tilting mutation has on the set of simple modules in the stable category.
  \end{abstract}

\maketitle
 
 \section{Introduction}
 
 Motivated by work of Okuyama and Rickard in modular representation theory, T. Aihara has recently introduced the notion of tilting mutation for symmetric algebras \cite{Aih1}.  Roughly speaking, these tilting mutations are defined as endomorphism rings of special tilting complexes of length $1$ that were first studied by Rickard \cite{DCSE} and Okuyama \cite{Oku} and have since proved quite useful in the construction of derived equivalences between blocks of finite groups as well as general symmetric algebras \cite{RFGT, Holm1, Holm2, TiltSym}.  These same tilting complexes have also been used by Vitoria \cite{Vit} and by Keller and Yang \cite{KellYang} to establish derived equivalences between the Jacobian algebras of certain pairs of quivers with potential which are linked by a mutation in the sense of Derksen, Weyman and Zelevinsky \cite{DWZ}.  Although these Jacobian algebras are often infinite-dimensional, this nice correspondence between the combinatorial mutation procedure of the quiver with potential and the homological mutation given by a derived equivalence warrants a further study of the combinatorics behind these derived equivalences in the finite-dimensional case.
 
The beginning of such a study is the first of two primary goals of the present article.  We aim to give a combinatorial description of tilting mutation for weakly symmetric algebras.  More precisely, given a weakly symmetric algebra $A$ presented by a quiver with relations and a vertex $i$ of the quiver at which there are no loops, we describe the quiver and relations of the endomorphism ring of a certain tilting complex associated to the vertex $i$.  This endomorphism ring, which we denote $\mu^+_i(A)$, is also a weakly symmetric algebra, and we say that it is obtained by mutating $A$ at the vertex $i$.  After reviewing some general facts about tilting mutation in Section 2, we describe the quiver of the mutated algebra in Section 3, and the relations for the mutated algebra in Section 4, both in terms of the quiver and relations for $A$.  We note that, in contrast to the mutation procedure for quivers with potential, in our setting the quiver for the mutated algebra typically depends on both the quiver and the relations of the original algebra.  Nevertheless, there remain some unsurprising similarities between the effects of tilting mutations and quiver mutations on the quiver of our algebra.  We see, for instance, reversal of arrows out of the vertex $i$ as well as new arrows corresponding to paths of length two through $i$.
 
 Secondly, we examine the effect of a tilting mutation inside the stable module category $\stmod A$.  Well-known work of Rickard has established that any two derived equivalent self-injective algebras are also stably equivalent  \cite{DCSE}.  However, little seems to be known about how the stable categories of two derived-equivalent self-injective algebras match up.  Studying this problem, we rediscovered a special case of an unpublished lemma of Okuyama \cite{Oku} describing the images in $\stmod A$ of the simple modules over the mutated algebra.   %As a consequence, we obtain a glimpse of how the stable categories of two often distinct weakly symmetric algebras, $A$ and $\mu^+_i(A)$, may non-trivially identify with one another in a large class of examples.  In fact, we have recently discovered that this result can be viewed as a special case of an unpublished lemma of Okuyama \cite{Oku}.  
 In his preprint, Okuyama applies this information to the problem of lifting a stable equivalence of Morita type to a derived equivalence in order to verify Brou\'{e}'s conjecture in several cases.  %The idea is that if one has an equivalence $\stmod B \rightarrow \stmod A$ of stable categories (of Morita type) which sends the simple $B$-modules to the images of the simple $\mu^+_i(A)$-modules, then a theorem of Linckelmann guarantees that $B$ and $\mu^+_i(A)$ are Morita equivalent, and hence that $B$ and $A$ are derived equivalent \cite{Lin}.
   In a similar vein, in Section 5 we apply our results and Okuyama's Lemma to resolve a question of Asashiba's concerning the derived Picard groups of a class of self-injective algebras of finite representation type \cite{Asa2}.  Specifically, for an algebra $A$ in this class, we show that a certain auto-equivalence of $\stmod A$, which was previously not known to be of Morita type, lifts to an auto-equivalence of the derived category $D^b(\rmod A)$ afforded by a tilting mutation.  Using Asashiba's arguments, the Main Theorem of \cite{Asa2} can now be strengthened to include all standard self-injective algebras of finite representation type as follows:
{\it  Any stable equivalence between standard representation-finite self-injective algebras lifts to a standard derived equivalence, and hence is of Morita type.}

Finally, motivated by the form Okuyama's Lemma takes in our setting (see Corollary 5.1), in Section 6 we abstract the effect of a tilting mutation in $\stmod A$ to a Hom-finite triangulated $k$-category $\T$ (with some additional hypotheses).  The main idea is to view the images of the simple $\mu^+_i(A)$-modules in $\stmod A$ as a mutation of the set of simple $A$-modules.  We develop this notion of mutation inside a triangulated category $\T$ for {\it maximal systems of orthogonal bricks}, which are sets of objects which homologically resemble the set of simple modules inside a stable category.  Our main result here is that the set of maximal systems of orthogonal bricks in $\T$ is closed under mutation.  In particular, we obtain a way to keep track of successive tilting mutations to an algebra $A$ inside the stable category of $A$ by successively mutating the set of simple $A$-modules.  Furthermore, such iterations of this mutation procedure will typically produce many nontrivial examples of maximal systems of orthogonal bricks in $\stmod A$, for a given algebra $A$.   In \cite{SMS} we will generalize this mutation procedure and compare it with related notions of mutation appearing in the work of Keller and Yang \cite{KellYang} and of Koenig and Yang \cite{KoYa}.
 
 Throughout this article, we work over a fixed algebraically closed field $k$.  We typically work with right modules and write morphisms on the left, composing them from right to left.  Likewise, paths in a quiver $Q$ will be composed from right to left, and we often identify them with morphisms between projective right modules over (a quotient of) the path algebra.  For an arrow $\alpha$ in a quiver, we shall write $s(\alpha)$ and $t(\alpha)$ for the source and target of $\alpha$ respectively.  Moreover, if $p$ and $q$ are paths in $Q$ we set 
$$ p/q = \left\{ \begin{array}{cc} p', & \mbox{if}\ p = p' q \\ 0, & \mbox{otherwise} \end{array} \right. \ \ \ \mbox{and}\ \ \ q \backslash p =  \left\{ \begin{array}{cc} p', & \mbox{if}\ p = q p' \\ 0, & \mbox{otherwise,} \end{array} \right.$$ and we extend these path-division operations linearly to $k$-linear combinations of paths $p$ and  $p'$ in the obvious way.  For an algebra $A$, we write $K(A)$ for the homotopy category of complexes of right $A$-modules.  We use complexes with differential of degree $1$ and define the degree-shifts of a complex $(X^{\bullet}, \delta)$ by $X[i]^p = X^{i+p}$ and $\delta[i] = (-1)^i\delta$ for $i \in \mathbb{Z}$.  We occasionally signal the degree-$0$ term of a complex by underlining it, and we identify $\rmod A$ with complexes concentrated in degree $0$.  We also frequently make use of a shorthand for matrices of morphisms, writing merely $[f_{ij}]$ in square brackets instead of a full matrix to denote a morphism $f : \oplus_{i \in I} X_i \rightarrow \oplus_{j \in J} Y_j$ where $f_{ij} : X_i \rightarrow Y_j$ for each $i, j$.  As is standard, we view elements of these direct sums as column vectors so that the morphism $f$ corresponds to left ``multiplication'' by the matrix $[f_{ij}]$.

\section{Tilting mutations and a lemma of Okuyama} 
\setcounter{equation}{0}

We assume that $A = k\Delta/I$ is a weakly symmetric $k$-algebra, presented as the path algebra of a quiver $\Delta$ modulo an admissible ideal $I$ of relations.  We let $J = k\Delta_{\geq 1}$ be the ideal of $k\Delta$ generated by the arrows, and write $J_A = J/I$ for the Jacobson radical of $A$.  We also write $\Delta_0 = \{1,2, \ldots, n\}$ and $\Delta_1$ for the vertices and arrows of $\Delta$ respectively.  For $U \subseteq \Delta_0$, we shall write $e_U = \sum_{i \in U} e_i$ for the sum of the corresponding primitive idempotents of $A$, and we shall write $P_U = e_U A$ for the corresponding projective $A$-module.  We also write $Q_U = (1-e_U)A$ so that $A_A \cong P_U \oplus Q_U$.  Letting $f_U : P_U \rightarrow L_U$ be a minimal left $\add(Q_U)$-approximation of $P_U$, it is not hard to see that $$T_U = [\ul{P_U} \stackrel{f_U}{\longrightarrow} L_U] \oplus Q_U[-1]$$ is a tilting complex concentrated in degrees $0$ and $1$.  Similarly, if $g_U : R_U \rightarrow P_U$ is a minimal right $\add(Q_U)$-approximation of $P_U$, then $${}_U T =  [R_U \stackrel{g_U}{\longrightarrow} \ul{P_U}] \oplus Q_U[1]$$ is a tilting complex concentrated in degrees $-1$ and $0$.

\begin{defin} Let $U \subseteq \Delta_0$.  The {\bf left (tilting) mutation of $A$ at $U$} is the algebra $$\mu^+_U(A) = \Endo_{K(A)}(T_U),$$ and the {\bf right (tilting) mutation of $A$ at $U$} is the algebra $$\mu^-_U(A) = \Endo_{K(A)}({}_U T).$$
\end{defin}

\noindent
{\it Remark.}  Our notation differs somewhat from Aihara's in \cite{Aih1}.  Namely, for a subset $U \subseteq \Delta_0$, Aihara defines a tilting complex $T(U)$, which is isomorphic to the complex ${}_{\bar{U}}T[1]$, where $\bar{U} = \Delta_0 \setminus U$.  It follows that the tilting mutation of $A$ at the vertex $i$, as defined by Aihara, coincides with the right tilting mutation $\mu^-_i(A)$ of $A$ at $i$ in the notation introduced above. \\

When we consider a sequence of such mutations, it is convenient to use the same indexing set for the vertices of the quivers of each mutated algebra. 
We employ the following convention.  The vertices of the quiver of $\mu^+_U(A)$ correspond to the indecomposable summands of $T_U$.  If $i \notin U$, then $e_i A$ (as a complex concentrated in degree $1$) is a summand of $T_U$ and we keep the label $i$ for the corresponding vertex of the new quiver.  If $i \in U$, then we use $i$ for the vertex corresponding to the summand $[\ul{P_i} \stackrel{f_i}{\longrightarrow} L_{U, i}]$ of $T_U$, where $f_i$ is a minimal left $\add(Q_U)$-approximation  (later we will instead write $i'$ for the vertex of the new quiver $\Delta'$ corresponding to $i \in \Delta_0$).

With these conventions we see that right and left mutations on the same subset of vertices yield inverse operations.

\begin{lemma} For any $U \subseteq \Delta_0$ we have $$\mu^-_U (\mu^+_U(A)) \cong A \cong  \mu^+_U (\mu^-_U(A)).$$ 
\end{lemma}

\noindent
{\it Proof.}  Let $G : D^b(A) \rightarrow D^b(B)$ be the equivalence induced by $T_U$, where $B = \mu^+_U(A) = \Endo_{K(A)}(T_U)$.  We write $P'_U$ and $Q'_U$ for the projective summands of $B$ given by $G([\ul{P_U} \rightarrow L_U])$ and $G(Q_U[-1])$ respectively.  The inverse equivalence then corresponds to the tilting complex $G(A)$ in $K(B)$.  To calculate $G(A)$, observe that we have the triangle
$$[\ul{P_U} \rightarrow L_U] \longrightarrow A \stackrel{\bnc{f_U \ \ 0}{0 \ \ 1}}{\longrightarrow} L_U \oplus Q_U \longrightarrow.$$ 
Applying $G$ yields a triangle in $K(B)$
$$P'_U \longrightarrow G(A) \stackrel{}{\longrightarrow} G(L_U \oplus Q_U) \longrightarrow$$ with $G(L_U \oplus Q_U) \in \add(Q'_U[1])$ since $L_U \oplus Q_U \in \add(Q_U) \subset \add(T_U[1])$.  It follows that $G(A)$ is concentrated in degrees $0$ and $-1$, with $P'_U$ in degree $0$ and with its degree-$(-1)$ term in $\add(Q'_U)$.  Since we know that $G(A)$ is a tilting complex, it follows that the component $G(L_U)[-1] \rightarrow P'_U$ of the connecting morphism in the above triangle must be a right $\add(Q'_U)$-approximation of $P'_U$.  Hence, $G(A)$ coincides with the tilting complex  ${}_U T'$ constructed over $B$ in the definition of $\mu^-_U(B)$.  The second isomorphism is proved symmetrically.  $\Box$ \\

For $U \subseteq \Delta_0$ fixed, we continue to let $B = \mu^+_U(A)$.  We let $F : D^-(\rmod B) \rightarrow D^-(\rmod A)$ be the equivalence constructed by Rickard in \cite{MTDC} that corresponds to the tilting complex $T_U \in K^b(\proj A)$, and we will denote by $\ul{F}$ the induced stable equivalence $\stmod B \rightarrow \stmod A$.  Following Rickard, to construct $F(X)$ for $X \in \rmod B$ one starts with a projective resolution $P^{\bullet}$ of $X$ and then uses the equivalence $\proj B \approx \add(T_U)$ to replace each term and each differential of this complex with the corresponding summands of $T \in K^b(\proj A)$ and the corresponding morphisms between them.  The result is nearly a double complex over $\proj A$, but the square of the differential is zero in only one direction, while in the other direction it is zero only up to homotopy.  Consequently, Rickard defines $F(X) \in K^-(\proj A)$ to be the total complex, but with differential given by a sum of maps of bi-degrees $(p,1-p)$ for $p \geq 0$.  Since $T$ is concentrated in degrees $0$ and $1$ the resulting complex $F(X)$ will be concentrated in degrees $\leq 1$.  Furthermore, Hu and Xi show that $F(X)$ has homology concentrated in degrees $0$ and $1$ and is isomorphic in $K^-(\rmod A)$ to a unique radical complex of the form $$0 \rightarrow Q^0 \rightarrow Q^1 \rightarrow 0$$ where $Q^1 \in \add(Q_U) = \add((1-e_U)A)$ \cite{HuXi3}.  One then sets $\ul{F}(X) := Q^0$ in $\stmod A$.

In his unpublished preprint \cite{Oku}, Okuyama gives the following description of the effect of the stable equivalence $\ul{F}$ on the simple $B$-modules, where we write $S_i$ (respectively, $S'_i$) for the simple $A$-module (resp. $B$-module) corresponding to $i \in \Delta_0$.

\vspace{2mm}
\noindent
{\bf Okuyama's Lemma.}  [Lemma 2.1' in \cite{Oku}]  \emph{\begin{itemize} \item[(1)] For $i \in U$, $\ul{F}(S'_i) \cong S_i$ in $\stmod A$.  
\item[(2)] For $j \notin U$, $\ul{F}(S'_j) \cong X$ in $\stmod A$, where $X$ is a submodule of $P_j$ such that \begin{itemize} \item[(i)] each composition factor of $X/\rad\ X$ is in the set $\{S_l\ |\ l \notin U\}$; and
\item[(ii)] each composition factor of $\rad\ \Omega^{-1}(X)$ is in the set $\{S_l\ |\ l \in U\}$.
\end{itemize}
\end{itemize}}

\vspace{2mm}
 While we do not repeat Okuyama's proof here, we explain in \cite{SMS} how an alternative proof can be derived from recent results of Koenig and Yang \cite{KoYa}.

\section{Arrows in the mutated quiver}
We continue to assume that $A = k\Delta/I$ is a weakly symmetric $k$-algebra.  For the remainder of this article, we focus on the case where $U = \{1\}$ consists of a single vertex $1 \in \Delta_0$ at which there are no loops.  To simplify our notation we shall henceforth replace the subscripts $U$ used above by $1$, or omit them altogether when there is no chance of confusion., and $1 \in \Delta_0 = \{1, \ldots, n\}$ is a vertex at which there are no loops.  Our first goal is to describe a quiver $\Delta'$ and an ideal $I'$ of relations in $k\Delta'$ giving a presentation of $B = \mu_1^+(A)$.

  As already remarked, the vertices of $\Delta'$ can be identified with the summands of $T$, and thus with the vertices of $\Delta$.  We shall write $\Delta'_0 = \{1', \ldots, n'\}$ where $1'$ corresponds to $T_1 := [\ul{P_1} \rightarrow L]$ and $i'$ corresponds to $T_i : = P_i[-1]$ for each $i \neq 1$.  Furthermore, notice that the former summand of $T$ has the form $P_1 \stackrel{[\gamma]}{\longrightarrow} \oplus_{s(\gamma)=1} P_{t(\gamma)}$.  We will also write $S'_i$ and $P'_i$ for the simple and indecomposable projective $B$-modules, respectively, associated to the vertex $i'$ of $\Delta'$.  
According to standard conventions, the arrows in $\Delta'$ from $i'$ to $j'$ should correspond to (a $k$-basis of) the irreducible maps in $\add(T)$ between the corresponding summands of $T$.  For simplicity, our initial description of the arrows in $\Delta'$ includes some maps which may turn out to be reducible, once the relations are taken into account.  Thus, the remainder of this section will focus on describing (the arrows of ) a quiver $\Delta'$, together with a surjective homomorphism $\Phi : k\Delta' \rightarrow \Endo_{K(A)}(T)$.  We will see that these arrows may arise in several possible ways, depending on the arrows of $\Delta$ as well as on certain relations in $I$.  

\begin{itemize}
\item[(A1)] Arrows $i' \rightarrow 1'$ for $i \neq 1$ are in one-to-one correspondence with the arrows $1 \rightarrow i$ in $\Delta$.  If $\alpha : 1 \rightarrow i$ in $\Delta$, we write $\alpha^* : i' \rightarrow 1'$ for the corresponding arrow in $\Delta'$, and we define $\Phi(\alpha^*)$ to be the irreducible map 
$$\xymatrix{0 \ar[r] \ar[d] & P_i \ar[d]^{1_{\alpha}} \\P_1 \ar[r]^-{[\gamma]} & \bigoplus_{s(\gamma)=1} P_{t(\gamma)},}$$ 
regarded as an endomorphism of $T$.

\vspace{2mm}
\item[(A2)]  Arrows $1' \rightarrow i'$ for $i \neq 1$ are in one-to-one correspondence with a fixed $k$-basis of \\ $(1-e_1)(I/(J(1-e_1)I+IJ))e_1$.  Given a nonzero map 
$$\xymatrix{ P_1 \ar[d] \ar[r]^-{[\gamma]} & \bigoplus_{s(\gamma)=1} P_{t(\gamma)} \ar[d]^{[f_{\gamma}]} \\ 0 \ar[r] & P_i}$$ 
we must first have $\sum_{\gamma} f_{\gamma} \gamma = 0$ in $A$.  That is, we must have a relation $r = \sum_{\gamma} f_{\gamma} \gamma \in I$, where the sum ranges over all arrows $\gamma$ starting in $1$ and $f_{\gamma}$ is a linear combination of paths from $t(\gamma)$ to $i$.  Since this map is assumed nonzero, we cannot have $r \in IJe_1$, and we may even assume each $f_{\gamma} \notin I$.  In order for this map to be irreducible, the relation $r$ must be left-minimal in the sense that it does not belong to $J (1-e_1)Ie_1$.  For any such left-minimal relation $r$, we shall write $r^*$ for the corresponding arrow $1' \rightarrow i'$ in $\Delta'$.  To obtain a complete set of such arrows $r^*$ for $\Delta'$, we need only let $r$ run through a $k$-basis of $(1-e_1)Ie_1/(1-e_1)(J(1-e_1)I+IJ)e_1 \cong (1-e_1)(I/(J(1-e_1)I+IJ))e_1$.  For each relation $r$ from this basis, we let $\Phi(r^*)$ be the irreducible map displayed above, where $f_{\gamma} = r/\gamma$.

\vspace{2mm}
\noindent
{\it Note:} When we speak of a basis of the quotient of an ideal in a path algebra $kQ$, such as $I/I'$ where $I' \subseteq I \subseteq kQ$, we really mean a set $\mathcal{B}$ of linear combinations of paths in $I \subseteq kQ$ whose images in $I/I'$ form a $k$-basis.  Moreover, we always assume that each $p \in \mathcal{B}$ is a linear combination of paths with the same source and target vertices in $Q$.

\vspace{2mm}
\item[(A3)] For $i, j \neq 1$, we have an arrow $i' \rightarrow j'$ in $\Delta'$ for each arrow $i \rightarrow j$ in $\Delta$.  If $\beta: i \rightarrow j$ is an arrow in $\Delta$, then we write $\beta' : i' \rightarrow j'$ for the corresponding arrow of $\Delta'$, and set $\Phi(\beta')$ equal to the map  $\beta[-1] : P_i[-1] \rightarrow P_j[-1]$ in $\add(T)$.  In certain cases this map may fail to be irreducible.  For example, in the relation $(r/\alpha)' - r^*\alpha^*$ introduced in (R2) of the next section, we would have $r/\alpha = \beta$ if $r := \beta \alpha$ is a relation in $I$.

\vspace{2mm}
\item[(A4)]  Also for $i,j \neq 1$, $\Delta'$ will contain a new arrow  $(\alpha \beta)' : i' \rightarrow j'$ for each pair of arrows $\beta: i \rightarrow 1$ and $\alpha : 1 \rightarrow j$ in $\Delta$ with $\alpha \beta \notin I$.  We set $\Phi((\alpha \beta)')$ equal to the map $\alpha \beta : e_iA[-1] \rightarrow e_jA[-1]$, which may be irreducible in $\add(T)$.  In certain cases, this map may still be reducible: for instance, if $I$ contains a relation of the form $\alpha \beta - p$ for some path $p$ that does not pass through $1$.  Alternatively, it is possible that the various irreducible maps obtained in this way are not linearly independent: for instance, if $I$ contains a relation of the form $\sum_i \alpha_i\beta_i$.  All of these possibilities will be handled by the relations $I'$ described in the next section.  (Thus our initial description of $I'$ is not necessarily an admissible ideal in the path algebra $k\Delta'$.)

\end{itemize}

This completes our description of the quiver $\Delta'$ and of the images in $\Endo_{K(A)}(T)$ of the arrows of $\Delta'$ under $\Phi$.  Naturally, we define $\Phi$ on the vertex $e'_i$ of $\Delta'$ to be the projection onto the corresponding summand of $T$.  It follows that $\Phi$ can be extended uniquely to an algebra homomorphism $k\Delta' \rightarrow \Endo_{K(A)}(T)$.

\begin{propos}  The homomorphism $\Phi : k\Delta' \rightarrow \Endo_{K(A)}(T) = B$ is surjective.
\end{propos}

\noindent
{\it Proof.}  It suffices to show that any (irreducible) radical map between indecomposable summands of $T$ can be factored through a linear combination of the maps $\Phi(a)$ for $a \in \Delta'_1$.  If $f : T_i \rightarrow T_1$ for $i \neq 1$, then clearly $f$ must factor through some linear combination of the maps $\Phi(\alpha^*)$ for $\alpha^*$ as in (A1).  Similarly, by construction, any $f : T_1 \rightarrow T_i$ for $i \neq 1$ as considered in (A2) will factor through some linear combination of the irreducible maps $\Phi(r^*)$ presented there.  Next, any radical map from $T_i$ to $T_j$ for $i, j \neq 1$ arises from a radical map $f : P_i \rightarrow P_j$ over $A$.  Clearly $f$ may be factored through the map 
$$\left( \begin{array}{c} [\beta] \\  \left[\alpha \beta\right] \end{array} \right) : P_i \rightarrow \bigoplus_{i \stackrel{\beta}{\rightarrow} l \neq 1} P_l \oplus \bigoplus_{i \stackrel{\beta}{\rightarrow} 1 \stackrel{\alpha}{\rightarrow} l \neq 1} P_l, $$
the components of which corrrespond via $\Phi$ to the arrows in $\Delta'$ starting at $i'$.

Finally, consider a radical morphism $f : P'_1 \rightarrow P'_1$:
$$\xymatrix{ P_1 \ar[d]^{f_0} \ar[r]^-{[\gamma]} & \bigoplus_{s(\gamma)=1} P_{t(\gamma)} \ar[d]^{f_1} \\ P_1 \ar[r]^-{[\gamma]} & \bigoplus_{s(\gamma)=1} P_{t(\gamma)}.}$$ 
Since $\Delta$ contains no loops at $1$, $f_0$ factors through the left approximation $[\gamma]$.  It follows that, up to homotopy, $f_0$ may be chosen to be $0$.  However, then the map $(0,f_1)$ will factor through $[\gamma^*]: \oplus_{s(\gamma)=1} P'_{t(\gamma)} \rightarrow P'_1$.  $\Box$ \\
%and hence it is not irreducible.

\section{Relations for the mutated algebra}

We now describe the ideal $I'$ of relations on the quiver $\Delta'$ so that $B = \mu^+(A) \cong k\Delta'/I'$.  These relations can emerge in various ways, and we divide them up based on their starting and ending vertices.    As mentioned already, the ideal $I'$ may fail to be contained in $k\Delta'_{\geq 2}$.  This means that some of the arrows in $\Delta'$ as described in the previous section may turn out to be redundant.  Nevertheless, it is convenient to include them here, as they lead to more uniform descriptions of $\Delta'$ and $I'$.

To describe the relations in $I'$, we will need a natural way of translating paths in $\Delta$ to paths in $\Delta'$.  Assume that $p$ is a path from $i$ to $j$ in $\Delta$ with $i, j \neq 1$.  We define the corresponding path $p'$ in $\Delta'$ inductively by the following rules:
\begin{itemize}
\item If $p = \beta$ is an arrow, then set $p' := \beta'$ as in (A3) in the previous section.
\item If $p = \alpha \beta$ for arrows $\beta : i \rightarrow 1$ and $\alpha : 1 \rightarrow j$ in $\Delta$, then set $p' := (\alpha \beta)'$ as in (A4) in the previous section.
\item If $p = p_1 p_2$ for a path $p_1$ that does not start at $1$, then set $p' := p_1'p_2'$.
\end{itemize}
A simple induction on the length of $p$ shows that $p'$ is well-defined.  Of course, we can extend this correspondence linearly to linear combinations of paths.  A less formal way of viewing this translation is by interpreting $p \in k\Delta$ as a map $P_i \rightarrow P_j$, which corresponds to a map $p': P_i[-1] \rightarrow P_j[-1]$ in $\add(T)$ and thus to a (not necessarily unique) linear combination of paths from $i'$ to $j'$ in $k\Delta'$.  However, by specifying the correspondence on the level of the path algebras, we can sidestep this issue of non-uniqueness.

%1) The last sentence rather describes an isomorphism $(1-e_1)A(1-e_1) \cong (1-e_1')B(1-e_1')$.  This means that the linear combination of paths in $k\Delta'$ is not necessarily unique. => Maybe this sentence should be removed.
%2) Does $\p mapsto p'$ induce a bijection $(1-e_1)k\Delta(1-e_1) \rightarrow (1-e_1')k\Delta'(1-e_1')$?

We now describe various relations in $k\Delta'$ and explain why each belongs to the kernel of $\Phi$.  By definition, $I'$ will be the ideal of $k\Delta'$ generated by these relations.  Afterwards, we will show that $I' = \ker \Phi$. \\

\begin{itemize}
\item[(R1)] Old relations $i' \rightarrow j'$ for $i, j \neq 1$:  If $\rho \in e_j Ie_i$, then the corresponding linear combination of paths $\rho'$ will be a relation from $i'$ to $j'$ in $I'$.  Notice that we may obtain minimal relations for $B$ in this way, even starting from non-minimal relations for $A$.  For instance, if $\rho : i \rightarrow 1$ is a minimal relation in $I$ and $\alpha : 1 \rightarrow j$, then $(\alpha \rho)' \in I'$ could be minimal.  In general, to obtain representatives of all the minimal relations of $I'$ that arise in this way, we need to consider a $k$-basis for $(1-e_1)[I / (J(1-e_1)I + I(1-e_1)J)](1-e_1)$.

\vspace{2mm}
\item[(R2)] Suppose $r$ is a (left minimal) relation in $e_iIe_1$ as in (A2) of the previous section.  We can decompose $r = \sum_{s(\alpha)=1}(r/\alpha)\alpha$ where $r/\alpha$ is a linear combination of paths from $t(\alpha)$ to $i$.  Then $I'$ contains the relation $$(r/\alpha)' - r^*\alpha^*$$ for each arrow $\alpha \in \Delta_1$ with source $1$.  To see that these relations are contained in $\ker \Phi$, observe that the composite
$$\xymatrix{P_{t(\alpha)}[-1] \ar[d]_{\Phi(\alpha^*)} & 0 \ar[r] \ar[d] & P_{t(\alpha)} \ar[d]^{1_{\alpha}} \\ T_1 \ar[d]_{\Phi(r^*)} & P_1 \ar[r]^-{[\gamma]} \ar[d] & \bigoplus_{s(\gamma)=1} P_{t(\gamma)} \ar[d]^{[r/\gamma]} \\ P_i[-1] & 0 \ar[r] & P_i}$$
clearly coincides with the map $\Phi((r/\alpha)') : P_{t(\alpha)}[-1] \rightarrow P_i[-1]$.  Notice that $(r/\alpha)'$ may be a single arrow or even zero.

\vspace{2mm}
\item[(R3)] For any arrow $\beta : i \rightarrow 1$, we have the following relation
$$\sum_{s(\alpha)=1} \alpha^*(\alpha \beta)' : i' \rightarrow 1' \in I'.$$
  To see that this is contained in $\ker \Phi$, notice that $\beta$ provides a homotopy between this sum and the zero map.
$$\xymatrix{P_{i}[-1] \ar[d]_{[\Phi((\alpha \beta)')]} & 0 \ar[r] \ar[d] & P_{i} \ar[d]^{[\alpha \beta]} \ar@{-->}[ddl]_(.3){\beta} \\ \bigoplus_{s(\alpha)=1} P_{t(\alpha)}[-1] \ar[d]_{[\Phi(\alpha^*)]} & 0 \ar[r] \ar[d] & \bigoplus_{s(\alpha)=1}P_{t(\alpha)} \ar[d]^{1} \\
T_1 & P_1 \ar[r]^-{[\alpha]} & \bigoplus_{s(\alpha)=1} P_{t(\alpha)} }$$

\vspace{2mm}
\item[(R4)] Suppose $\rho: 1 \rightarrow 1$ is a minimal relation in $I$, so that $\alpha \rho : 1 \rightarrow i$ is a relation inducing an arrow $(\alpha \rho)^*$ (as in (A2) of the previous section) for any arrow $\alpha : 1 \rightarrow i$.  Then we have $$\sum_{s(\alpha)=1} \alpha^*(\alpha \rho)^* : 1' \rightarrow 1'  \in I'.$$  In fact the corresponding morphism $T_1 \rightarrow T_1$ is null-homotopic via the map $[\rho/\gamma]$ where $\rho = \sum_{s(\gamma)=1} (\rho/\gamma) \gamma$.
$$\xymatrix{T_1 \ar[d]_{[\Phi((\alpha \rho)^*)]} & P_1 \ar[r]^-{[\gamma]} \ar[d] & \bigoplus_{s(\gamma)=1} P_{t(\gamma)} \ar[d]^{[\alpha \rho / \gamma]} \ar@{-->}[ddl]_(.3){[\rho/\gamma]} \\ \bigoplus_{s(\alpha)=1} P_{t(\alpha)}[-1] \ar[d]_{[\Phi(\alpha^*)]} & 0 \ar[r] \ar[d] & \bigoplus_{s(\alpha)=1} P_{t(\alpha)} \ar[d]^{1} \\ T_1 & P_1 \ar[r]^-{[\alpha]} & \bigoplus_{s(\alpha)=1} P_{t(\alpha)}}$$

\vspace{2mm}
\item[(R5)] The relations from $1' \rightarrow i'$ with $i \neq 1$ are the most difficult to describe explicitly.  We can identify them by way of the following lemma.
\begin{lemma} A linear combination $\rho$ of paths from $1'$ to $i'$ is contained in $\ker \Phi$ if and only if $\rho \alpha^* \in \ker \Phi$ for all arrows $\alpha \in \Delta_1$ with source $1$ (if and only if the map $P_{t(\alpha)} \rightarrow P_i$ that corresponds to $\rho \alpha^*$ is zero in $A$ for all such $\alpha$).
\end{lemma}

\noindent
{\it Proof.}  The forward direction is clear.  Thus assume that $\rho\alpha^* \in \ker \Phi$ for all arrows $\alpha : 1 \rightarrow i$ in $\Delta$.  Since the arrows of the form $\alpha^*$ are the only arrows with target $1'$ in $\Delta'$, we see that $\rho J_B = 0$ and thus $\rho \in \soc\ P'_1$.  As $B$ is weakly symmetric, we must have $\soc\ P'_1 \cong S'_1$, which forces $i = 1$ (a contradiction) or $\rho = 0$ in $B$. $\Box$ \\

\noindent
Thus we add to $I'$ all linear combinations $\rho$ of paths from $1'$ to $i'$ for which we already have $\rho \alpha^* \in I'$ for all arrows $\alpha \in \Delta_1$ with source $1$. \\

\noindent
{\it Remark.}  Let $\rho : 1' \rightarrow i'$ be a minimal relation arising as above.  Observe that by using the relations from (R2), we may replace any $r^*\alpha^*$ occuring in $\rho$ by $(r/\alpha)'$.  We may thus choose a basis of $(1-e'_1)(I'/(J'I'+I'J'))e'_1$ consisting of linear combinations of paths not containing any arrow $r^*$ other than as the initial arrow.
\end{itemize}

\begin{propos} We have $I' = \ker \Phi$.  Hence $\Phi$ induces an isomorphism $k\Delta'/I' \cong \Endo_{K(A)}(T) = B$.
\end{propos}

\noindent
{\it Proof.}  Let $\rho \in \ker \Phi$.  Without loss of generality we may assume $\rho \in e'_i k\Delta' e'_j$ for idempotents $e'_i$ and $e'_j$ corresponding to vertices $i'$ and $j'$ of $\Delta'$, respectively.  %Furthermore, it suffices to treat the case where $\rho$ does not belong to $J'(\ker \Phi) + (\ker \Phi) J'$.  
%By induction on the minimal length of a path with nonzero coefficient in $\rho$, we may further assume that $\rho$ does not belong to $J'(\ker \Phi) + (\ker \Phi) J'$.  

First, if $i, j \neq 1$, then notice that (R2) allows us to add elements of $I'$ to $\rho$ to eliminate any occurrences of arrows of the form $r^*$ and $\alpha^*$ as in (A2) and (A1).  We may thus assume that $\rho = p'$ for some linear combination of paths $p$ from $i$ to $j$ in $k\Delta$.  Then the fact that $\Phi(\rho) = 0$ is equivalent to $p \in I$.  Thus $\rho \in I'$ according to (R1).

If $i=1$ and $j \neq 1$ (cf. (R5)), we clearly have $\rho \alpha^* \in \ker \Phi$ for all arrows $\alpha$ with source $1$. Thus $\rho \alpha^* \in I'$ for all such $\alpha$ by the above argument  and $\rho \in I'$ by (R5).

If $i \neq 1$ and $j=1$ (cf. (R3)), we can decompose $\rho = \sum_{s(\alpha)=1} \alpha^* (\alpha^* \backslash \rho)$, and (R2) again allows us to assume that $\alpha^* \backslash \rho$ is of the form $p_{\alpha}'$ for a linear combination of paths $p_{\alpha} \in k\Delta$.  For $\Phi(\rho)$ to be null-homotopic, the map $[p_{\alpha}] = [\Phi(\alpha^* \backslash \rho)[1]] : P_i \rightarrow \oplus_{s(\alpha)=1} P_{t(\alpha)}$ must factor through $[\alpha] : P_1 \rightarrow \oplus_{s(\alpha)=1} P_{t(\alpha)}$, meaning that we can write $p_{\alpha} = \alpha p$ for a single linear combination of paths $p$ from $i$ to $1$.  We then have $$\rho = \sum_{s(\alpha)=1} \alpha^* \sum_{t(\beta)=1}(\alpha \beta)' (\beta \backslash p)'  \in I'$$ by (R3).

The case where $i=j=1$ (cf. (R4)) is treated similarly.  We begin by observing that we may assume $\rho \notin J'(\ker \Phi)$ and that $\rho \notin (\ker \Phi)J'$ by the two previous cases.  Now, for each $\alpha, \gamma \in \Delta_1$ with source $1$, the above argument applied to $\rho \gamma^*$ shows that $\alpha^* \backslash \rho \gamma^* = (\alpha p_{\gamma})'$ for some linear combination of paths $p_{\gamma} \in k\Delta$ from $t(\gamma)$ to $1$.  Thus the factorization $\Phi(\rho) = [\Phi(\alpha^*)][\Phi(\alpha^* \backslash \rho)]$ corresponds to the following maps in $\add(T)$, which is null-homotopic via $[p_{\gamma}]$.

$$\xymatrix{T_1 \ar[d]_{[\Phi(\alpha^* \backslash \rho)]} & P_1 \ar[r]^-{[\gamma]} \ar[d] & \bigoplus_{s(\gamma)=1} P_{t(\gamma)} \ar[d]^{[\alpha p_{\gamma}]} \ar@{-->}[ddl]_(.3){[p_{\gamma}]} \\ \bigoplus_{s(\alpha)=1} P_{t(\alpha)}[-1] \ar[d]_{[\Phi(\alpha^*)]} & 0 \ar[r] \ar[d] & \bigoplus_{s(\alpha)=1} P_{t(\alpha)} \ar[d]^{1} \\ T_1 & P_1 \ar[r]^-{[\alpha]} & \bigoplus_{s(\alpha)=1} P_{t(\alpha)}}$$
Since the vertical composite on the left is zero, we must have $q := \sum_{s(\gamma)=1} p_{\gamma} \gamma \in e_1 I e_1$.  Next we check that $q$ is a minimal relation so that the construction in (R4) can be applied.  Our assumption that $\rho \notin J'(\ker \Phi)$ implies that $[\Phi(\alpha^* \backslash \rho)] \neq 0$ and hence that $\alpha p_{\gamma} \notin I$ for some $\alpha$ and $\gamma$.   Consequently, $p_{\gamma} \notin I$ for some $\gamma$, and $q \notin IJ$.   On the other hand, if $q \in JI$, then we have $q = \sum_{t(\beta)=1} \beta (\beta \backslash q)$ with each $\beta \backslash q \in I$.  In the above diagram, we can factor the degree-$(-1)$ component of $\Phi(\rho)$ as $[\alpha p_{\gamma}] = [\alpha \beta] [ \beta \backslash p_{\gamma}] = [\alpha] [\beta] [\beta \backslash p_{\gamma}]$.   As $[\beta] : \oplus_{t(\beta)=1} P_{s(\beta)} \rightarrow P_1$ gives a homotopy between the map $(0\ \  [\alpha \beta]) : \oplus_{t(\beta)=1} T_{s(\beta)} \rightarrow T_1$ and the zero map, we would have  $\rho \in (\ker \Phi)J'$, contradicting our earlier assumption.   
 
Finally, for each arrow $\gamma$ with source $1$, we have $$\alpha^* \backslash \rho \gamma^* = (\alpha p_{\gamma})' = (\alpha q / \gamma)' \equiv (\alpha q)^* \gamma^* \ (\mbox{mod}\ I')$$ by (R2).   Thus, by (R5), we have $\alpha^* \backslash \rho \equiv (\alpha q)^*\ (\mbox{mod}\ I')$, and hence $$\rho = \sum_{s(\alpha)=1} \alpha^* (\alpha^* \backslash \rho) \equiv \sum_{s(\alpha)=1} \alpha^* (\alpha q)^* \equiv 0 \ (\mbox{mod}\ I')$$
according to (R4).  This completes the proof that $\ker \Phi \subseteq I'$.  $\Box$ \\

\vspace{2mm}
As an application of the above results we can describe the beginnings of (not necessarily minimal) projective resolutions of the simple right $B$-modules.  %We will use these resolutions in the next section to compute the images of these simples under the stable equivalence induced by the tilting complex $T$.  
We shall apply the general description of a projective resolution of a simple module $S_i = e_i(kQ/I)$ over an algebra $\Lambda = kQ/I$, which begins
 $$\bigoplus_{p \in e_i(I/(IJ + JI))} e_{s(p)}\Lambda  \stackrel{[\alpha \backslash p]}{\longrightarrow} \bigoplus_{\alpha \in Q_1,\ t(\alpha)=i} e_{s(\alpha)}\Lambda \stackrel{[\alpha]}{\longrightarrow} e_i \Lambda \longrightarrow S_i \rightarrow 0,$$
 where the first direct sum is indexed by a $k$-basis of $e_i(I/(IJ+JI))$.  This resolution is actually minimal if $I$ is an admissible ideal (i.e., if the presentation of $\Lambda$ by quiver with relations is minimal), but that is not always the case for the quiver and relations $(Q',I')$ described above.

Thus a projective resolution of $S'_1$ begins 

\begin{eqnarray}
 \bigoplus_{p \in e_1(I/(IJ+JI))e_1} P'_1 \oplus \bigoplus_{t(\beta)=1} P'_{s(\beta)} \stackrel{\left( [(\alpha p)^*]\ [(\alpha \beta)'] \right)}{\longrightarrow} \bigoplus_{s(\alpha)=1} P'_{t(\alpha)} \stackrel{[\alpha^*]}{\longrightarrow} P'_1 \longrightarrow S'_1 \rightarrow 0,
\end{eqnarray}
where again the first direct sum is indexed by a $k$-basis of the specified set.  The first map between projectives is determined by the arrows of type (A1) from Section 3, and the second map is determined by the minimal relations described in (R3) and (R4) above.  Similarly, a projective resolution for $S'_i$ with $i \neq 1$ begins

\begin{eqnarray}
 \bigoplus_{\rho \in e'_iI'e'_1} P'_1 \oplus \bigoplus_{q \in e_iIe_1} \left[ \bigoplus_{s(\delta)=1} P'_{t(\delta)} \right] \oplus \bigoplus_{p \in e_iI(1-e_1)} P'_{s(p)} & \stackrel{\varphi}{\longrightarrow} & \nonumber \\ \bigoplus_{r \in e_iIe_1} P'_1 \oplus \bigoplus_{j \stackrel{\gamma}{\rightarrow} 1 \stackrel{\alpha}{\rightarrow} i} P'_{s(\gamma)} \oplus \bigoplus_{t(\beta)=i} P'_{s(\beta)} & \stackrel{\left( [r^*]\ [(\alpha \gamma)']\ [\beta'] \right)}{\longrightarrow} & P'_i \longrightarrow S'_i \rightarrow 0,
\end{eqnarray}
where the map $\varphi$ is given by
 $$\varphi =  \left( \begin{array}{ccc} 0 & -[\delta_{r,q} \delta^*] & 0\\  \left[(\alpha \gamma)' \backslash \rho \right] &  \left[(\alpha \gamma \backslash q/ \delta)'\right] & \left[(\alpha \gamma \backslash p)'\right] \\ \left[\beta' \backslash \rho \right] & \left[ (\beta \backslash q/ \delta)' \right] & \left[ (\beta \backslash p)' \right] \end{array} \right) .$$
Furthermore, in the summations $\rho$ runs through a $k$-basis of $e'_i( I'/(J'I'+I'J')) e'_1$ consisting of relations as in (R5), $q$ and $r$ run through the same $k$-basis of $e_i(I/(J(1-e_1)I +IJ))e_1$ as in (A2) of the last section (hence we may write $\delta_{r,q}$ for the Kronecker delta), and $p$ runs through a $k$-basis of $e_i[I / (J(1-e_1)I + I(1-e_1)J)](1-e_1)$ as in (R1) above.  Observe that the $0$ in the upper-left entry of $\varphi$ is a consequence of the choice of basis described in the Remark following Lemma 4.1.  In particular, for any $\rho$ belonging to such a basis of $e'_i(I'/(J'I'+I'J'))e'_1$ and any $r$, we have $r^* \backslash \rho = 0$.

\section{Examples and Applications}

Now that we have described the quiver and relations of the mutated algebra $B = \mu_1^+(A)$, which is derived equivalent to the original algebra $A$, our goal in the remainder of this paper is to apply Okuyama's Lemma to study the stable equivalence between $A$ and $B$ that is induced by the tilting complex $T$.  As above, we write $S_1, \ldots, S_n$ and $S'_1, \ldots, S'_n$ for the simple modules (up to isomorphism) over $A$ and $B$ respectively, and we let $\ul{F} : \stmod B \rightarrow \stmod A$ be the induced stable equivalence.  In our context, we obtain the following specialization of Okuyama's Lemma.

\begin{coro} With notation as above, we have $\ul{F}(S'_1) \cong S_1$ and for each $i \neq 1$, $\ul{F}(S'_i) \cong e_i J_A (1-e_1)A$, which can also be described as the largest submodule of $\rad\ P_i$ without a factor of $S_1$ in its top Loewy layer.
\end{coro}

For our first example, we consider a standard self-injective algebra of finite representation type which is of type  $(\mathbb{D}_{3m}, 1/3, 1)$ for $m \geq 2$.   In fact, by Asashiba's classification theorem \cite{DECSA}, any such algebra is derived equivalent to the algebra $A$ presented by the quiver
$$
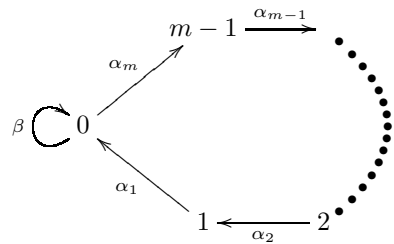
\xymatrix{ & m-1 \ar@{<-}[dl]_{\alpha_m} &   \ar@{<-}[l]_-{\alpha_{m-1}} \\ 0 \ar@(dl, ul)^{\beta} \ar@{<-}[dr]_{\alpha_1} \\ & 1 \ar@{<-}[r]_{\alpha_2}  & 2 \ar @{{}{*}} @/_2pc/[uu]}$$
and relations (i) $\alpha_1 \cdots \alpha_m = \beta^2$; (ii) $\overbrace{\alpha_i \cdots \alpha_m \beta \alpha _1 \cdots  \alpha_i}^{m+2} = 0$ for all $i \in \{1, \ldots, m\} = \mathbb{Z}/\gen{m}$; and (iii) $\alpha_m \alpha_1 = 0$ \cite{Asa2}. 

Our interest in this specific example stems from the following problem raised by Asashiba in \cite{Asa2}, where it is shown that nearly all equivalences between the stable categories of standard self-injective algebras of finite representation type can be lifted to standard derived equivalences.  (By a {\it standard} derived equivalence we mean a derived equivalence induced by tensoring with a two-sided tilting complex as in \cite{DEDF}).  Asashiba shows that there is essentially one stable auto-equivalence for the algebras of type $(\mathbb{D}_{3m},s/3,1)$ with $3\nmid s$ for which this problem is left unresolved.  More specifically, the stable AR-quiver of such an algebra $\Lambda$ is isomorphic to $\mathbb{ZD}_{3m}/\gen{\tau^{(2m-1)s}}$, and the stable module category $\stmod \Lambda$ is (the $k$-variety generated by) the mesh category $k(\mathbb{ZD}_{3m}/\gen{\tau^{(2m-1)s}})$ since $\Lambda$ is standard.  Hence, corresponding to the order-$2$ graph automorphism of $\mathbb{D}_{3m}$, we have an automorphism $H$ of the stable AR-quiver of $\Lambda$ and therefore of the category $\stmod \Lambda$ as well, which fixes most indecomposables, but swaps each pair of indecomposables corresponding to the leafs of a fixed sectional $\mathbb{D}_{3m}$ subquiver or one of its $\tau$-translates.  By Asashiba's description of the stable picard group of $\Lambda$ (i.e., the group of all auto-equivalences of $\stmod \Lambda$, modulo automorphisms), to show that this auto-equivalence is induced by an auto-equivalence of $D^b(\rmod \Lambda)$ it suffices to show that at least one stable equivalence which does not induce a power of $\tau$ on the stable AR-quiver can be lifted to a derived equivalence.  For the algebra $A$ given above ($s=1$), we will show that such a stable equivalence $\ul{F}$ is induced by a tilting mutation, and hence by a derived equivalence.  We will then indicate how to extend this result to the $\mathbb{Z}/s$-Galois covering of $A$ which gives a representative algebra of type $(\mathbb{D}_{3m},s/3,1)$.  Finally, by Corollary 3.5 in \cite{DEDF}, we can conclude that there is a standard derived equivalence having the same effect on isomorphism classes of objects, and hence inducing a stable equivalence that agrees with $\ul{F}$ on isomorphism classes of objects.

For the algebra $A$, we now check what happens when we perform a tilting mutation at the vertex $1$.  The quiver of the resulting endomorphism algebra contains arrows $\beta'$ and $\alpha_i'$ for $3 \leq i \leq m$.  We also have new arrows $(\alpha_1 \alpha_2)'$, $\alpha_1^*$, $(\alpha_m \alpha_1)^*$ and $r^*$ where $r$ can be taken to be %either $\alpha_1 \cdots \alpha_m \beta \alpha_1$ or 
$\beta^2\alpha_1$.

$$
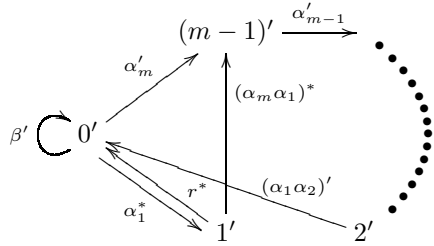
\xymatrix{ & (m-1)' \ar@{<-}[dl]_{\alpha'_m} &   \ar@{<-}[l]_-{\alpha'_{m-1}} \\ 0' \ar@(dl, ul)^{\beta'}  \ar@{<-}[drr]^(.7){(\alpha_1\alpha_2)'} \ar[dr]<-1ex>_{\alpha_1^*}\\ & 1' \ar[ul]_(.3){r^*}  \ar[uu]_(.7){(\alpha_m\alpha_1)^*} & 2' \ar @{{}{*}} @/_2pc/[uu]}$$

For the relations, we first consider those induced by the original relations of $A$.  From (i) we get $$(\alpha_1\alpha_2)' \alpha'_3 \cdots \alpha'_m= (\beta')^2;$$ from (ii) we get $$\alpha'_i \cdots \beta' (\alpha_1 \alpha_2)' \cdots \alpha'_i = 0 \ \mbox{for}\ 3 \leq i \leq m$$ as well as $$(\alpha_1\alpha_2)' \cdots \alpha'_m \beta' (\alpha_1 \alpha_2)' = 0 $$ which is induced by the corresponding non-minimal relation from $2$ to $0$.  
From (R2) applied to the relations $r$ and $\alpha_m \alpha_1$ we obtain 
\begin{eqnarray*} (\beta')^2 & =  & r^*\alpha_1^* \\
 \alpha'_m & = & (\alpha_m \alpha_1)^*\alpha_1^*.
\end{eqnarray*}
From (R3) we obtain $$\alpha_1^*(\alpha_1 \alpha_2)' = 0.$$
We do not obtain any relations from (R4), but (R5) yields the relation
$$r^* =  (\alpha_1 \alpha_2)' \cdots \alpha'_{m-1} (\alpha_m \alpha_1)^*$$ since we have $r^*\alpha_1^* = (\beta')^2 = (\alpha_1\alpha_2)' \cdots \alpha'_m = (\alpha_1\alpha_2)' \cdots \alpha'_{m-1} (\alpha_m\alpha_1)^*\alpha_1^*$.  Thus the mutated algebra $B$ can be presented by the quiver 
$$
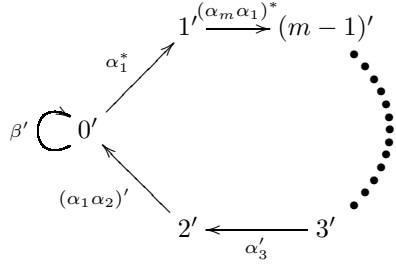
\xymatrix{ & 1'  \ar@{<-}[dl]_{\alpha_1^*}  \ar[r]^-{(\alpha_m\alpha_1)^*} &   (m-1)' \\ 0'  \ar@(dl, ul)^{\beta'} \ar@{<-}[dr]_{(\alpha_1\alpha_2)'} \\ & 2' \ar@{<-}[r]_{\alpha'_3}  & 3' \ar @{{}{*}} @/_2pc/[uu]}$$
which is isomorphic to the quiver of $A$.  Moreover, the relations above correspond precisely to the relations of $A$ under such an isomorphism of quivers.  Thus we have $\mu_1^+(A) \cong A$ in this example.
 
In light of this isomorphism, we may consider the induced derived equivalence $F$ as an autoequivalence of $D^b(A)$, which induces the autoequivalence $\und{F}$ of $\stmod A$.  Since this isomorphism identifies $0'$ with $0$, $1'$ with $m-1$ and $i'$ with $i-1$ for $2 \leq i \leq m-1$, Corollary 5.1 shows that $\und{F}$ has the following effect on the simple $A$-modules 
\begin{eqnarray*} \und{F}(S_{m-1}) & \cong & S_1 \cong \Omega^3(S_{m-1}), \\ \und{F}(S_0) & \cong & e_0J_A(1-e_{1})A \cong \Omega^3(e_0A/\beta \alpha_1 A), \\ \und{F}(S_i) & \cong & e_{i+1}J_A(1-e_1)A \cong \Omega(S_{i+1})  \cong \Omega^3(S_i) \ \mbox{for}\ 1 \leq i \leq m-2. \end{eqnarray*}

We claim that no power of the syzygy functor $\Omega$ (or, equivalently, of $\tau$) can have the same effect on the simples in $\stmod A$.   This is due to the fact that $\Omega$ coincides with $\tau^{m}$ on objects, while $S_0$ and $e_0A/\beta \alpha_1 A$ are not in the same $\tau$-orbit.  The latter assertion can be seen by the fact that $S_0$ and $e_0A/\beta \alpha_1 A$ occur as summands (shown below as the first two) of the middle term in the same almost split sequence:

$$\xymatrixcolsep{.8pc} \xymatrixrowsep{1pc} 
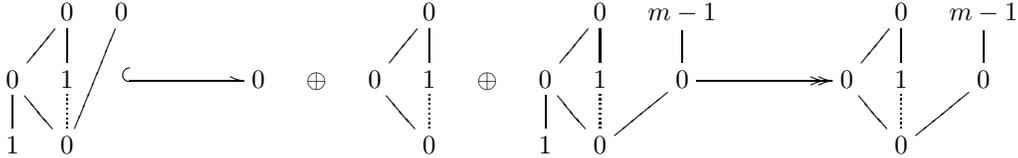
\xymatrix{ & 0 \ar@{-}[dl] \ar@{-}[d]  & 0 \ar@{-}[ddl]   &  & &  &  & & 0 \ar@{-}[dl] \ar@{-}[d] &  & & 0 \ar@{-}[d] \ar@{-}[dl] & m-1 \ar@{-}[d]    & & & & 0 \ar@{-}[dl] \ar@{-}[d] & m-1 \ar@{-}[d]  \\ 
 0 \ar@{-}[dr] \ar@{-}[d] & 1 \ar@{{}{.}}[d] & \ar@^{(->}[rrr] & & &   0  &     \oplus   & 0 \ar@{-}[dr] & 1 \ar@{{}{.}}[d] &   \oplus  & 0 \ar@{-}[d] \ar@{-}[dr] & 1 \ar@{{}{.}}[d] & 0 \ar@{-}[dl]  \ar@{->>}[rrr] & & & 0 \ar@{-}[dr] & 1 \ar@{{}{.}}[d] & 0 \ar@{-}[dl] \\
 1                                  & 0               &                      & & &     &            & & 0                                          &      & 1 & 0 & & & & & 0 }$$
 Since the AR-quiver of $A$ is isomorphic to $\mathbb{ZD}_{3m}/\gen{\tau^{2m-1}}$, the indecomposable summands of the middle term of any almost split sequence must lie in distinct $\tau$-orbits. 

We now consider the $\mathbb{Z}/s$-grading on $A$ where $\deg (\beta) = 1, \deg (\alpha_1) = 2$ and $\deg(\alpha_i) = 0$ for $2 \leq i \leq m$, and we let $\tilde{A}$ be the corresponding Galois covering of $A$.  It is known that $\tilde{A}$ is a standard self-injective algebra of finite representation type $(\mathbb{D}_{3m},s/3,1)$, the category $\rmod \tilde{A}$ is equivalent to the category $\gr A$ of finitely-generated $\mathbb{Z}/s$-graded $A$-modules, and the forgetful functor $\gr A \rightarrow \rmod A$ is dense.  With respect to this grading of $A$, the tilting complex $T = [e_1 A \stackrel{\alpha_1}{\rightarrow} e_0 A] \oplus (1-e_1)A[-1]$ is gradable (with $e_1 A$ generated in degree $-2$ and all other projectives generated in degree $0$), and thus the mutated algebra $\mu^+_1(A) = \Endo_{K(A)}(T)$ inherits a $\mathbb{Z}/s$-grading.  Moreover, one easily checks that the isomorphism $\mu^+_1(A) \cong A$ respects this grading, as the map $(\alpha_1 \alpha_2)'$ comes from an endomorphism of degree $2$, $\beta'$ comes from an endomorphism of degree $1$ and all the other arrows in the mutated quiver are derived from endomorphisms of degree $0$.  Now, according to Theorem 4.4 of \cite{GreHap}, the tilting complex $T$ can be lifted to $\tilde{A}$, thereby yielding a derived auto-equivalence $\tilde{F}$ of $\tilde{A}$ that lifts the derived auto-equivalence $F$ constructed above.  Finally, if $\tilde{F}$ induced a functor on the stable category $\stmod \tilde{A} \approx \stgr A$ that agreed with a power of the syzygy functor on isomorphism classes of objects, the same would have to be true of $F$ on $\stmod A$ since the forgetful functor $\stgr A \rightarrow \stmod A$ is dense.  However, we have already shown that this is not the case.

\vspace{5mm}
As a second example, we consider the algebra $A$ given by the quiver
$$\xymatrixrowsep{3pc} 
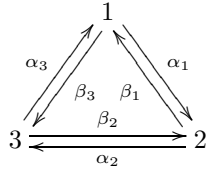
\xymatrix{ & 1 \ar[dr]<.5ex>^{\alpha_1} \ar[dl]<.5ex>^{\beta_3} \\ 3 \ar[rr]<.5ex>^{\beta_2} \ar[ur]<.5ex>^{\alpha_3} & & 2 \ar[ll]<.5ex>^{\alpha_2} \ar[ul]<.5ex>^{\beta_1}}$$
and the relations $ \beta_i \alpha_i = 0 = \alpha_i \beta_i$ and $\alpha_{i+2}\alpha_{i+1}\alpha_i = \beta_i \beta_{i+1} \beta_{i+2}$ for each $i\  (\mbox{mod}\ 3)$.
Performing a tilting mutation at the vertex $1$, we obtain  arrows as follows:
\begin{itemize}
\item[(A1)] Corresponding to the arrows $\alpha_1$ and $\beta_3$, we have $\alpha_1^* : 2' \rightarrow 1'$ and $\beta_3^* : 3' \rightarrow 1'$.
\item[(A2)] Corresponding to the relations $\beta_3 \beta_1 \alpha_1 = 0$, $\alpha_1 \alpha_3 \beta_3 = 0$, $\alpha_1 \alpha_3 \alpha_2 \alpha_1 = 0$ and $\beta_3 \beta_1 \beta_2 \beta_3 = 0$ we obtain arrows $(\beta_3 \beta_1 \alpha_1)^*: 1' \rightarrow 3' ;\  (\alpha_1 \alpha_3 \beta_3)^*: 1' \rightarrow 2'; \ (\alpha_1 \alpha_3 \alpha_2 \alpha_1)^* : 1' \rightarrow 2'$ and $(\beta_3 \beta_1 \beta_2 \beta_3)^* : 1' \rightarrow 3'$.
\item[(A3)] Corresponding to the arrows $\alpha_2$ and $\beta_2$, we have $\alpha'_2 : 2' \rightarrow 3'$ and $\beta'_2 : 3' \rightarrow 2'$.
\item[(A4)] Corresponding to the nonzero paths of length $2$ passing through $1$, we have the arrows $(\alpha_1 \alpha_3)' : 3' \rightarrow 2'$ and $(\beta_3 \beta_1)' : 2' \rightarrow 3'$.
\end{itemize}

Following Section 4, we work out the relations to be:
\begin{itemize}
\item[(R1)] From the relations in $A$ between vertices $2$ and $3$, we obtain $\beta_2' \alpha'_2 = \alpha_2' \beta'_2 = 0, \ (\alpha_1 \alpha_3)' \alpha'_2 = \beta'_2 (\beta_3 \beta_1)'$ and $\alpha'_2 (\alpha_1 \alpha_3)' = (\beta_3 \beta_1)' \beta'_2$.

\item[(R2)] The relations of the form $(r/\alpha)' = r^* \alpha^*$ where $r$ runs through the relations used in (A2) above are $(\beta_3 \beta_1)' = (\beta_3 \beta_1 \alpha_1)^* \alpha_1^*, 0 = (\beta_3 \beta_1 \alpha_1)^* \beta_3^*; (\alpha_1 \alpha_3)' = (\alpha_1 \alpha_3 \beta_3)^* \beta_3^*, 0 =  (\alpha_1 \alpha_3 \beta_3)^* \alpha_1^*; (\alpha_1 \alpha_3)' \alpha'_2 = (\alpha_1 \alpha_3 \alpha_2 \alpha_1)^* \alpha_1^*;  0 = (\alpha_1 \alpha_3 \alpha_2 \alpha_1)^* \beta_3^*; (\beta_3 \beta_1)' \beta'_2 = (\beta_3 \beta_1 \beta_2 \beta_3)^*\beta_3^*; 0 = (\beta_3 \beta_1 \beta_2 \beta_3)^* \alpha_1^*$.  In particular, we see that the arrows $(\beta_3 \beta_1)'$ and $(\alpha_1 \alpha_3)'$ may be eliminated from the quiver.

\item[(R3)] For $\beta = \alpha_3, \beta_1$, we get the relations $\alpha_1^*(\alpha_1 \alpha_3)' = 0$ and $\beta_3^* (\beta_3 \beta_1)' = 0$ respectively.

\item[(R4)] Taking $\rho = \alpha_3 \beta_3$ and $\rho = \beta_1 \alpha_1$ yields $\alpha_1^* (\alpha_1 \alpha_3 \beta_3)^* = 0$ and $\beta_3^* (\beta_3 \beta_1 \alpha_1)^* = 0$ respectively.  Taking $\rho = \alpha_3 \alpha_2 \alpha_1 - \beta_1 \beta_2 \beta_3$ yields $\alpha_1^* (\alpha_1 \alpha_3 \alpha_2 \alpha_1)^* - \beta_3^*(\beta_3 \beta_1 \beta_2 \beta_3)^* = 0$.

\item[(R5)] Combining relations from (R2), (R1) and (R2) again, we have $(\alpha_1 \alpha_3 \alpha_2 \alpha_1)^*\alpha_1^* = (\alpha_1 \alpha_3)' \alpha_2' = \beta'_2 (\beta_3 \beta_1)' = \beta'_2 (\beta_3 \beta_1 \alpha_1)^* \alpha_1^*$ and $(\alpha_1 \alpha_3 \alpha_2 \alpha_1)^* \beta_3^* = 0 = \beta'_2(\beta_3 \beta_1 \alpha_1)^* \beta_3^*$ by (R2).  Thus we have $(\alpha_1 \alpha_3 \alpha_2 \alpha_1)^* = \beta'_2(\beta_3 \beta_1 \alpha_1)^*$ using Lemma 4.1.  Similarly, one checks that $(\beta_3 \beta_1 \beta_2 \beta_3)^* = \alpha'_2 (\alpha_1 \alpha_3 \beta_3)^*$.

\end{itemize}

\vspace{2mm}
Thus, after eliminating the redundant arrows, the quiver of the mutated algebra becomes
$$\xymatrixrowsep{4pc} \xymatrixcolsep{3pc} 
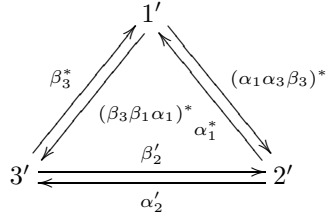
\xymatrix{ & 1' \ar[dr]<.5ex>^{(\alpha_1 \alpha_3 \beta_3)^*} \ar[dl]<.5ex>^{(\beta_3 \beta_1 \alpha_1)^*} \\ 3' \ar[rr]<.5ex>^{\beta'_2} \ar[ur]<.5ex>^{\beta_3^*} & & 2' \ar[ll]<.5ex>^{\alpha'_2} \ar[ul]<.5ex>^(.4){\alpha_1^*}}$$
and the relations show that this algebra is in fact isomorphic to the original algebra $A$.  By Corollary 5.1, we have $$ \und{F}(S'_1) \cong S_1;\ \  \und{F}(S'_2) \cong \xymatrixrowsep{.5pc} \xymatrixcolsep{1pc} \xymatrix{3 \ar@{-}[d] & 3 \ar@{-}[ddl] \\ 1 \ar@{-}[d] \\ 2}; \ \ \und{F}(S'_3) \cong \xymatrixrowsep{.5pc} \xymatrixcolsep{1pc} \xymatrix{2 \ar@{-}[d] & 2 \ar@{-}[ddl] \\ 1 \ar@{-}[d] \\ 3}.$$

%What happens if we mutate at a pair of vertices?  Is the mutated algebra still isomorphic to this one?  If so, are there any other basic algebras derived equivalent to this algebra?

\section{Mutations of maximal systems of orthogonal bricks}

In Sections 2 and 3 we have described the quiver and relations for the algebra obtained via a tilting mutation at a single vertex without a loop.  Even in this restricted setting, it would be interesting to understand the changes to the quiver and relations resulting from performing successive tilting mutations.  While this problem may be too general to admit a nice answer, we can still apply Okuyama's Lemma to keep track of these successive mutations internally in $\stmod A$ by following the images of the simple modules.  In fact, this leads us to a way of mutating sets of objects that behave like simple modules in the stable category $\stmod A$.

For starters, suppose that $A$ is a weakly symmetric algebra with simples $S_1, \ldots, S_n$ such that $\Ext^1_A(S_1,S_1)=0$, and let $B = \mu^+_1(A)$.   If $\und{F} : \stmod B \rightarrow \stmod A$ is the induced stable equivalence, Corollary 5.1 shows that the images of the simple $B$-modules under $\und{F}$ are the modules $S^+_i$ defined by $S^+_1 = S_1$ and for $i \neq 1$ by the short exact sequences $$\ses{S^+_i}{\Omega S_i}{S_1^{n_i}}{}{f_i}$$ where $f_i$ is a minimal left $\add(S_1)$-approximation.  In fact, we can replace these short exact sequence by distinguished triangles $S^+_i \longrightarrow S_i[-1] \stackrel{\und{f_i}}{\longrightarrow} S_1^{n_i} \rightarrow$ in the stable category $\stmod A$, where $-[d]$ denotes the $d^{th}$ power of the suspension functor $\Omega^{-1}$.  This observation motivates us to define a general mutation procedure on sets of objects resembling simples in certain triangulated categories.

We thus let $\T$ be a Hom-finite triangulated Krull-Schmidt $k$-category with suspension denoted by $-[1]$.  In order for $\T$ to resemble the stable module category of a weakly symmetric algebra, we assume additionally that $$ \T(X,Y) \cong  \T(Y,X[-1])$$ for all $X, Y \in \T$.  These isomorphisms are not required to be natural, and as these Hom-spaces are only $k$-vector spaces, this is equivalent to asking that these two Hom-spaces always have equal dimensions.  Such an equality holds, for instance, if $\T$ is $-1$-Calabi-Yau.  %In fact, $\stmod A$ is $-1$-Calabi-Yau precisely when $A$ is symmetric (and not semi-simple).
  Following Pogorza\l y \cite{Pog}, we say that a (finite) set $\s = \{S_1, \ldots, S_n\}$ of objects of $\T$ is a {\it maximal system of orthogonal bricks} if the following hold:
\begin{enumerate} 
\item $\T(S_i,S_i) \cong k$ for all $i$.
\item $\T(S_i,S_j) = 0$ for all $i \neq j$.
\item For every nonzero $X \in \T$ there exists a $j$ such that $\T(X, S_j) \neq 0$.
\item $S_i[2] \not \cong S_i$ for all $i$.
\end{enumerate}
Note that condition (3) is equivalent to (3') For every $X \in \T$ there exists a $j$ such that $\T(S_j, X) \neq 0$ in light of the identity $\dim_k \T(X,Y) = \dim_k \T(Y, X[-1])$.   We have also substituted condition (4) for the requirement (from Pogorza\l y's original definition) that $\tau S_i \not \cong S_i$ for all $i$, as these two are equivalent for the set of simples in the stable category of a weakly symmetric algebra.

\begin{defin}  Suppose $\s = \{S_1, \ldots, S_n\}$ is maximal system of orthogonal bricks in $\T$ such that $\T(S_i,S_i[1]) = 0$.  We define the {\bf left mutation of $\s$ at $S_i$} to be $\mu^+_i(\s) := \{S^+_1, \ldots, S^+_n\}$ where $S^+_i = S_i$ and for $j \neq i$, $S^+_j$ is defined via  a distinguished triangle $$S^+_j \rightarrow S_j[-1] \stackrel{f_j}{\rightarrow} S_i^{n_j} \rightarrow$$
where $f_j$ is a minimal left $\add(S_i)$-approximation.  In particular, since $\T(S_i,S_i) \cong k$, the map $f_j$ has the form $(g_1, \ldots, g_{n_j})^T : S_j[-1] \rightarrow S_i^{n_j}$ where $\{g_1, \ldots, g_{n_j}\}$ is a $k$-basis for $\T(S_j[-1],S_i)$.
\end{defin}

\noindent
{\it Remark.}  An equivalent mutation of simples appears in the work of Keller and Yang on Jacobian algebras of quivers with potential  \cite{KellYang}.  They construct an equivalence of derived categories (induced by an Okuyama complex) that has the same effect on the simples once it is composed with the suspension functor.  We thank Bernhard Keller for bringing our attention to this result.\\

We first show that left mutation preserves the defining properties of maximal systems of orthogonal bricks.  

\begin{therm} Let $\s = \{S_1, \ldots, S_n\}$ be a maximal system of orthogonal bricks in $\T$, and suppose that $\T(S_i, S_i[1]) = 0$.  Then the set $\mu^+_i(\s)$ is again a maximal system of orthogonal bricks in $\T$.
\end{therm}  

\noindent
{\it Proof.}  Applying $\T(-,S_i)$ to the triangle $S^+_j \rightarrow S_j[-1] \stackrel{f_j}{\rightarrow} S_i^{n_j} \rightarrow$ yields an exact sequence $$\T(S_i^{n_j},S_i) \rightarrow \T(S_j[-1],S_i) \stackrel{0}{\rightarrow} \T(S^+_j,S_i) \rightarrow \T(S_i^{n_j}[-1],S_i),$$
where the first map is an epimorphism by definition of $f_j$ and the last term is $0$ since $S_i$ is assumed to have no self-extensions.   Thus $\T(S^+_j,S_i)=0$.  Similarly applying $\T(S_i,-)$ to the same triangle yields the exact sequence.
$$\xymatrix{\T(S_i,S^+_j[-1]) \ar[r] \ar[d]^{\cong} &  \T(S_i, S_j[-2]) \ar[r] \ar[d]^{\cong} & \T(S_i,S_i^{n_j}[-1]) \ar[r] \ar[d]^{\cong} & \T(S_i, S^+_j) \ar[r] & \T(S_i,S_j[-1]) \ar[d]^{\cong} \\ \T(S^+_j,S_i) \ar@{=}[d] & \T(S_j[-1],S_i) \ar@{=}[d] & \T(S_i^{n_j},S_i) \ar@{=}[d] & & \T(S_j,S_i) \ar@{=}[d] \\ 0 & k^{n_j} & k^{n_j} & & 0}$$
We thus see that $\T(S_i, S^+_j) = 0$.
Now suppose that $j, l \neq i$, and consider the exact diagram obtained by applying $\T(S^+_j,-)$ and $\T(-,S_l[-1])$ to the triangles defining $S^+_l$ and $S^+_j$ respectively.
$$\xymatrix{ & & \T(S_i^{n_j}[-1], S_l[-1]) \\ \T(S^+_j, S_i^{n_l}[-1]) \ar[r]^0 &\T(S^+_j, S^+_l) \ar[r] & \T(S^+_j, S_l[-1]) \ar[r]^0 \ar[u]^0 & \T(S^+_j, S_i^{n_l}) \\ & & \T(S_j[-1], S_l[-1]) \ar[u] \\ & & \T(S_i^{n_j}, S_l[-1]) \ar[u]^0 
}$$
To get the zero maps, we have used that $\T(S^+_j, S_i^{n_l}[-1]) \cong \T(S_i^{n_l},S^+_j)=0$, $\T(S_i,S_l)=0$, $\T(S^+_j,S_i)=0$ and $\T(S_i^{n_j},S_l[-1]) \cong \T(S_l,S_i^{n_j}) = 0$.  Thus the remaining nonzero maps now yield isomorphisms $$\T(S^+_j, S^+_l) \cong \T(S^+_j, S_l[-1]) \cong \T(S_j[-1],S_l[-1]) \cong \T(S_j,S_l) \cong \left\{ \begin{array}{cc} k & j = l \\ 0 & j \neq l. \end{array} \right.$$

Next we check that condition (3) holds for $\mu^+_i(\s)$.  Assume that $\T(X,S^+_j) = 0$ for all $j$ for some nonzero $X \in \T$.  In particular, $\T(X,S_i) =0$ since $S^+_i = S_i$.  Thus for all $j \neq i$, applying $\T(X,-)$ to the defining triangle for $S^+_j$ shows that $\T(X,S_j[-1]) = 0$.  Consequently, $\T(S_j, X) = 0$ for all $j \neq i$, which forces $\T(S_i,X) \neq 0$ by the dual of condition (3) for the maximal system of orthogonal bricks $\s$.  Now let $g: S_i^{r} \rightarrow X$ be a minimal right $\add(S_i)$-approximation of $X$, where $r = \dim_k \T(S_i,X)>0$, and define $Y$ as the cone of $g$ in a distinguished triangle $$S_i^r \stackrel{g}{\longrightarrow} X \longrightarrow Y \rightarrow.$$
Since $\T(S_i[1],S^+_j) \cong \T(S^+_j, S_i) = 0$ for all $j \neq i$, we see that $\T(Y,S^+_j) = 0$ for all $j \neq i$.  Next, we apply $\T(S_i,-)$ to the above triangle to get 
$$\xymatrix{\T(S_i,X[-1]) \ar[r] \ar@{=}[d] &  \T(S_i, Y[-1]) \ar[r]  & \T(S_i,S_i^{r}) \ar[r]^{\T(S_i,g)} \ar[d]^{\cong} & \T(S_i, X) \ar[d]^{\cong} \ar[r] & \T(S_i,Y) \ar@{=}[d] \ar[r]& \T(S_i, S_i^r[1]) \ar@{=}[d] \\ 0&  & k^r &  k^r & 0 & 0 }$$
The first term above is zero since it is isomorphic to $\T(X,S_i)=0$, the last term is zero since $\T(S_i,S_i[1])=0$.  Furthermore, $\T(S_i,g)$ must be an isomorphism since it is an epimorphism by definition, while its domain and co-domain have the same $k$-dimension.  It follows that $\T(Y,S_i) \cong \T(S_i,Y[-1]) = 0$.  Now $Y$ satisfies $\T(Y,S^+_j) = 0$ for all $j$, and hence as we argued above for $X$, we must have $\T(S_i,Y) \neq 0$ or else $Y=0$.   We have already seen that the first possibility does not occur, so suppose $Y=0$.  Then, by definition of $Y$, $X \cong S_1^r$, which contradicts $\T(X,S_1)=0$.  

Finally, for condition (4), observe that $\T(S_i[2],S_i) \cong \T(S_i,S_i[1])=0$, which implies that $S_i[2] \not \cong S_i$.  If $j \neq i$ and $n_j >0$, then the exact sequence $\T(S^+_j[2],S_i) \rightarrow \T( S_i^{n_j}[1],S_i) \rightarrow \T(S_j,S_i)$, obtained from applying $\T(-,S_i)$ to the triangle defining $S^+_j$, yields $\T(S^+_j[2],S_i) \neq 0$, as $\T(S_j,S_i) = 0$ while $\T(S_i^{n_j}[1],S_i) \cong \T(S_i,S_i^{n_j}) \neq 0$.  Thus $S^+_j[2] \not \cong S^+_j$.  Alternatively, if $n_j=0$, then $S^+_j = S_j[-1]$ and $S^+_j[2] \not \cong S^+_j$ follows from $S_j[2] \not \cong S_j$.  $\Box$ \\

\noindent
{\it Remarks.}  (1) Recently, Koenig and Liu have defined {\it simple-minded systems} of objects (in stable categories) \cite{KoeLiu}, and the definition carries over easily to a triangulated category $\T$ as above.  Furthermore, they note that any simple minded system in this setting is also a maximal system of orthogonal bricks; although the converse is not clear.  It turns out that the set of simple-minded systems in $\T$ is also closed under the mutations defined above, but the argument requires a closer study of torsion pairs in $\T$ and hence will be presented in a sequel to this article \cite{SMS}.

(2) It is shown in \cite{CX1} that stable equivalences of Morita type (between self-injective algebras) preserve {\it rigid} modules, i.e., modules with open $\GL_d(k)$-orbits in the affine variety $\Mod^A_d$ parametrizing $d$-dimensional $A$-modules.  As simple modules are always rigid, one can construct many further examples of rigid modules by successively mutating the simple $A$-modules in the above manner.\\

We can also define right mutations on the sets of maximal systems of orthogonal bricks in $\T$.  Suppose $\s = \{S_1, \ldots, S_n\}$ is a maximal system of orthogonal bricks in $\T$ with $\T(S_i, S_i[1])=0$.  Then the {\bf right mutation of $\s$ at $S_i$} is defined as the set $\mu^-_i(\s) = \{S^-_1, \ldots, S^-_n\}$ where $S^-_i = S_i$ and for $j \neq i$, $S^-_j$ is defined via a distinguished triangle $$S_i^{m_j}  \stackrel{g_j}{\rightarrow} S_j[1] \rightarrow S^-_j \rightarrow$$ where $g_j$ is a minimal right $\add(S_i)$-approximation.  As expected, left and right mutation at $S_i$ turn out to be inverse operations.  To see this, it suffices to observe that in the triangle defining $S^+_j \in \mu^+_i(\s)$, the connecting morphism $S_i^{n_j} \rightarrow S^+_j[1]$ is a right $\add(S_i)$-approximation (since its cone is $S_j$ and $\T(S_i,S_j)=0$).  In particular, $S^{+-}_j \in \mu^-_i(\mu^+_i(\s))$ is defined as the cone of this connecting morphism, which is just $S_j$.  Hence we have the following.

\begin{propos} For any maximal system of orthogonal bricks $\s$ in $\T$ such that $\T(S_i,S_i[1])=0$, we have $\mu^+_i(\mu^-_i(\s)) = \s = \mu^-_i(\mu^+_i(\s))$.
\end{propos}

We now explain how these mutations keep track of the images of the simple modules under the stable equivalences induced by successive tilting mutations.  More specifically, suppose that $F : D^b(A') \rightarrow D^b(A)$ is an equivalence of triangulated categories, inducing an equivalence $\und{F}: \stmod A' \rightarrow \stmod A$, and $\s = \{\und{F}(S'_1), \ldots, \und{F}(S'_n)\}$ where the $S'_i$ are the simple $A'$-modules.  Setting  $A'' = \mu^+_i(A')$, we get a stable equivalence $\und{G} : \stmod A'' \rightarrow \stmod A'$, and assuming $\Ext^1_{A'}(S'_i,S'_i)=0$,  Corollary 5.1 implies that the image of the set of simple $A''$-modules $\s'' = \{S''_1, \ldots, S''_n\}$ under $\und{G}$ is $\mu^+_i(S'_1, \ldots, S'_n)$.  Thus $\und{F}(\und{G}(\s'')) = \und{F}(\mu^+_i(S'_1, \ldots, S'_n)) = \mu^+_i(\s)$, since $\und{F}$ preserves triangles and approximations.  In fact, by general arguments involving the preservation of exact sequences by a stable equivalence (see for example \cite{IPUSE}), the same is true if $\und{F}$ is replaced by any stable equivalence between algebras without nodes, even one that is not assumed to be an equivalence of triangulated categories.

\begin{coro} Let $A$ be a weakly symmetric algebra, and let $\M_{st}$ (respectively, $\M_{der}$) be the set of all maximal systems of orthogonal bricks in $\stmod A$ which correspond to the simple $B$-modules via a stable equivalence $\stmod B \rightarrow \stmod A$ (respectively, via a stable equivalence induced by a derived equivalence $\D^b(B) \rightarrow D^b(A)$) for some algebra $B$.  Then $\M_{st}$ (resp. $\M_{der}$) is closed under left and right mutation.
\end{coro}


\begin{thebibliography}{99}


\bibitem{Aih1} T. Aihara.  \emph{Mutating Brauer trees.}  Preprint (2011) arXiv:1009.3210v3.

%\bibitem{Aih2} T. Aihara.  \emph{Silting mutation for self-injective algebras.}  Preprint (2010) arXiv:1012.3265v1.

%\bibitem{AihIya} T. Aihara and O. Iyama.  \emph{Silting mutation in triangulated categories.}  Preprint (2011) arXiv:1009.3370v3.

%\bibitem{Amiot} C. Amiot.  \emph{On the structure of triangulated categories with finitely many indecomposables.}   Bull. Soc. Math. France  135  (2007),  no. 3, 435--474.

\bibitem{DECSA} H. Asashiba.  \emph{The derived equivalence classification of representation-finite selfinjective algebras.}  J. Algebra 214 (1999), no. 1, 182--221.

\bibitem{Asa2} H. Asashiba.  \emph{On a lift of an individual stable equivalence to a standard derived equivalence for representation-finite self-injective algebras.}  Algebr. Represent. Theory 6 (2003), no. 4, 427--447.

%\bibitem{Cover} H. Asashiba.  \emph{Covering functors, skew group categories and derived equivalences.}  Preprint (2008) arXiv:0807.4706v2.

%\bibitem{SEAA} M. Auslander and I. Reiten.  \emph{Stable equivalence of Artin algebras.}  Proc. of the conf. on orders, group rings and related topics (Ohio State Univ., Columbus, Ohio, 1972), 8-71, Lecture Notes in Math., 353. Springer, Berlin, 1973.

%\bibitem{DPA} J. Bia\l kowski, K. Erdmann and A. Skowro\'{n}ski.  \emph{Deformed preprojective algebras of generalized Dynkin type.}  Trans. Amer. Math. Soc. 359 (2007), no. 6, 2625-2650.

%\bibitem{GCY3} R. Bocklandt.  \emph{Graded Calabi-Yau algebras of dimension 3.} (Appendix by M. Van den Bergh) J. Pure Appl. Algebra 212 (2008), no. 1, 14--32.

%\bibitem{PAAK} S. Brenner, M. C. R. Butler and A. D. King.  \emph{Periodic algebras which are almost Koszul.}  Algebras and Representation Theory  5 (2002), 331-367.

%\bibitem{BMRRT} A. B. Buan, R. Marsh, M. Reineke, I. Reiten, and G. Todorov.  \emph{Tilting theory and cluster combinatorics.}  Adv. Math. 204 (2006), no. 2, 572--618.

%\bibitem{Buch} R. O. Buchweitz.  \emph{Finite representation type and periodic Hochschild (co-)homology.}  Trends in the representation theory of finite-dimensional algebras (Seattle, WA, 1997),  81--109, Contemp. Math., 229, Amer. Math. Soc., Providence, RI, 1998. 

%\bibitem{BIKR} I. Burban, O. Iyama, B. Keller and I. Reiten.  \emph{Cluster tilting for one-dimensional hypersurface singularities.}   Adv. Math.  217  (2008),  no. 6, 2443--2484.

\bibitem{RFGT} J. Chuang and J. Rickard.  \emph{Representations of finite groups and tilting.} Handbook of tilting theory, 359--391, London Math. Soc. Lecture Note Ser., 332, Cambridge Univ. Press, Cambridge, 2007.

\bibitem{DWZ} H. Derksen, J. Weyman and A. Zelevinsky.  \emph{Quivers with potentials and their representations. I. Mutations.} Selecta Math. (N.S.) 14 (2008), no. 1, 59--119.

%\bibitem{Per} A. Dugas.  \emph{Periodic resolutions and self-injective algebras of finite type.}  Preprint (2008) arXiv:0808.1311v2.


\bibitem{TiltSym} A. Dugas. \emph{A construction of derived equivalent pairs of symmetric algebras.}  Preprint (2010) arXiv:1005.5152v2 [math.RT].

\bibitem{CX1} A. Dugas.  \emph{On periodicity in bounded projective resolutions.}  Preprint (2012) arXiv:1203.2408v1 [math.RT].


\bibitem{SMS} A.  Dugas.  \emph{Torsion pairs and simple-minded systems in triangulated categories.}  Preprint (2012) arXiv:1207.7338 [math.RT].

%\bibitem{Erd} K. Erdmann. \emph{Blocks of tame representation type and related algebras.} Lecture Notes in Mathematics, 1428. Springer-Verlag, Berlin, 1990.

%\bibitem{ErdSko}  K. Erdmann and A. Skowro\'nski.  \emph{The stable Calabi-Yau dimensionof tame symmetric algebras.}  J. Math. Soc. Japan 58 (2006), 97-128.

%\bibitem{ErdSko2} K. Erdmann and A. Skowro\'nski.  \emph{Periodic algebras.}  Trends in Representation Theory and Related Topics.  European Math. Soc., Zurich, 2008.


%\bibitem{Green} E. L. Green.  \emph{Graphs with relations, coverings and group-graded algebras.}  Trans. Amer. Math. Soc.  279  (1983),  no. 1, 297--310.

\bibitem{GreHap} E. L. Green and D. Happel.  \emph{Gradings and derived categories.}  Algebr. Represent. Theor. 14 (2011) 497--513.

%\bibitem{DCFDA} D. Happel.  \emph{On the derived category of a finite-dimensional algebra.}  Comment. Math. Helv. 62 (1987), no. 3, 339-389.

%\bibitem{TCRTA} D. Happel.  \emph{Triangulated categories in the representation theory of finite-dimensional algebras.} London Mathematical Society Lecture Note Series, 119. Cambridge University Press, Cambridge, 1988.

%\bibitem{HCFDA} D. Happel.  \emph{Hochschild cohomology of finite-dimensional algebras.}  S\'eminaire d'Alg\`{e}bre Paul Dubreil et Marie-Paul Malliavin, 39\`{e}me Ann\'ee (Paris, 1987/1988),  108--126, Lecture Notes in Math., 1404, Springer, Berlin, 1989.

\bibitem{Holm1} T. Holm.  \emph{Derived equivalent tame blocks.} J. Algebra 194 (1997) 178--200.

\bibitem{Holm2} T. Holm. \emph{Derived equivalence classification of algebras of dihedral, semidihedral, and quaternion type.}  J. Algebra 211 (1999) 159--205.

\bibitem{HoKa} M. Hoshino and Y.  Kato.  \emph{Tilting complexes defined by idempotents.}  Comm. Algebra 30 (2002), no. 1, 83--100. 

%\bibitem{HuKX} W. Hu, S. Koenig and C. Xi.  \emph{Derived equivalences from cohomological approximations, and mutations of perforated Yoneda algebras.} In preparation.


%\bibitem{HuXi1} W. Hu and C. Xi.  \emph{Almost $\mathcal{D}$-split sequences and derived equivalences.}  Preprint (2008) arXiv:0810.4757v1 [math.RT].

\bibitem{HuXi3} W. Hu and C. Xi.  \emph{Derived equivalences and stable equivalences of Morita type, I.}  Nagoya Math. J. 200 (2010) 107--152.
%Preprint (2008) arXiv:0810.4761v2 [math.RT].

%\bibitem{HuXi2} W. Hu and C. Xi.  \emph{Derived equivalences for $\Phi$-Auslander-Yoneda algebras.} Preprint (2010) arXiv:0912.0647v2.

%\bibitem{Iyama1} O. Iyama.  \emph{Higher-dimensional Auslander-Reiten theory on maximal orthogonal subcategories.} Adv. Math. 210 (2007), no. 1, 22--50.

%\bibitem{Iyama2} O. Iyama.  \emph{Auslander correspondence.} Adv. Math. 210 (2007), no. 1, 51--82.

%\bibitem{IR} O. Iyama and I. Reiten.  \emph{Fomin-Zelevinsky mutation and tilting modules over Calabi-Yau algebras.}  Amer. J. Math.  130  (2008),  no. 4, 1087--1149.


%\bibitem{TOC} B. Keller.  \emph{On triangulated orbit categories.}  Documenta Math. 10 (2005), 551-581.

%\bibitem{DCT} B. Keller. \emph{Derived categories and tilting.} Handbook of tilting theory, 49--104, London Math. Soc. Lecture Note Ser., 332, Cambridge Univ. Press, Cambridge, 2007.

%\bibitem{CYTC} B. Keller.  \emph{Calabi-Yau triangulated categories.}  Trends in Representation Theory of Algebras and Related Topics, European Math. Soc., Zurich, 2008.

\bibitem{KellYang} B. Keller and D. Yang.  \emph{Derived equivalences from mutations of quivers with potential.}  Adv. Math. 226 (2011) 2118--2168.

\bibitem{KoeLiu} S. Koenig and Y. Liu.  \emph{Simple-minded systems in stable module categories.}  Preprint (2010) arXiv:1009.1427v1.

\bibitem{KoYa} S. Koenig and D. Yang.  \emph{Silting objects, simple-minded collections, $t$-structures and co-$t$-structures for finite-dimensional algebras.}  Preprint (2012) arXiv:1203.5657v2.

\bibitem{Lin} M. Linckelmann.  \emph{Stable equivalences of Morita type for self-injective algebras and $p$-groups.}  Math. Z. 223 (1996), no. 1, 87--100. 


\bibitem{IPUSE} R. Mart\'{i}nez Villa.  \emph{Properties that are left invariant under stable equivalence.}  Comm. Algebra 18 (1990), no. 12, 4141--4169.

%\bibitem{DPic} J. Miyachi and A. Yekutieli.  \emph{Derived Picard groups of finite-dimensional hereditary algebras.}  Composito Mathematica 129 (2001), 341-368.

%\bibitem{NATC} A. Neeman. \emph{Some new axioms for triangulated categories.} J. Algebra  139  (1991),  no. 1, 221--255. 

\bibitem{Oku} T. Okuyama.  \emph{Some examples of derived equivalent blocks of finite groups.}  Preprint (1998).

\bibitem{Pog} Z. Pogorza\l y.   \emph{Algebras stably equivalent to self-injective special biserial algebras.}   Comm. in Algebra 22 (1994), no. 4, 1127--1160.

\bibitem{MTDC} J. Rickard.  \emph{Morita theory for derived categories.}  J. London Math. Soc. 39 (1989), no. 2, 436--456.

\bibitem{DCSE} J. Rickard.  \emph{Derived categories and stable equivalence.}  J. Pure Appl. Algebra 61 (1989), no. 3, 303--317.

\bibitem{DEDF} J. Rickard.  \emph{Derived equivalences as derived functors.}  J. London Math. Soc. 43 (1991), no. 2, 37--48.

%\bibitem{EDCSA} J. Rickard.  \emph{Equivalences of derived categories for symmetric algebras.} J. Algebra 257 (2002), no. 2, 460--481.

%\bibitem{ADK} Ch. Riedtmann.  \emph{Algebren, Darstellungsk\"ocher,  Ueberlagerungen und zur\"uck.}  Comment. Math. Helvetici 55 (1980), 199--224.

\bibitem{Vit} J. Vitoria.  \emph{Mutations vs. Seiberg duality.}  J. Algebra 321  (2009),  no. 3, 816--828. 

\end{thebibliography}
\end{document}